\documentclass[12pt]{article}\def\today{December 7, 2000}
\usepackage{amssymb}
\usepackage{amsthm,amsmath}
\newtheorem{thm}{Theorem}[section]
\newtheorem{cor}[thm]{Corollary} 
\newtheorem{lem}[thm]{Lemma}

\newtheorem{definition}[thm]{Definition}
\newenvironment{de}{\begin{definition}\rm}{\end{definition}}
\newtheorem{example}[thm]{Example}
\newenvironment{exmp}{\begin{example}\rm}{\end{example}}
\newtheorem{remark}[thm]{Remark}

\newcommand{\Section}[1]{\section{#1}\setcounter{equation}{0}}

\newcommand{\zwei}[2]{\left[
    \begin{array}{c} #1 \\ #2
    \end{array} \right]} \newcommand{\vier}[4]{\left[
    \begin{array}{ccc} #1 &\;& #2 \\ #3 &\;& #4 \end{array}
  \right]}

\newcommand{\C}{{\mathbb C}}
\newcommand{\R}{{\mathbb R}}
\newcommand{\pp}{{\mathbb P}}
\newcommand{\PN}{{\mathbb P}^N}

\newcommand{\eqr}[1]{~\mbox{$(${\rm \ref{#1}}$)$}}

\title{Output Feedback Invariants\thanks{The main results of this
    paper were announced without proof at the 36th IEEE
    Conference on Decision and Control in San Diego, 1997.
    (see~\cite{ra97p}).}} \date{\today}
\author{M. S. Ravi\\
  {\small Department of Mathematics}\\
  {\small East Carolina University}\\
  {\small Greenville, NC 27858}\\
  {\small ravim@mail.ecu.edu} \and Joachim
  Rosenthal\thanks{Supported
    in part by NSF grant DMS-96-10389.}\\
  {\small Department of Mathematics}\\
  {\small University of Notre Dame}\\
  {\small Notre Dame, Indiana 46556-5683}\\
  {\small Rosenthal.1@nd.edu} \and
Uwe Helmke\thanks{ Research partially supported
          by DFG grant 436 RUS 113/275/4-1.}\\
  {\small Department of Mathematics}\\
  {\small University of W\"urzburg}\\
  {\small 97074 W\"urzburg, Germany}\\
  {\small helmke@mathematik.uni-wuerzburg.de}}

\begin{document}
\maketitle
\begin{abstract} 
  The paper is concerned with the problem of determining a
  complete set of invariants for output feedback. Using tools
  from geometric invariant theory it is shown that there exists a
  quasi-projective variety whose points parameterize the output
  feedback orbits in a unique way. If the McMillan degree $n\geq
  mp$, the product of number of inputs and number of outputs,
  then it is shown that in the closure of every feedback orbit
  there is exactly one nondegenerate system.
\end{abstract}\newpage

\Section{Introduction}        \label{Sec:Int}

Consider a time invariant linear (strictly proper) system
\begin{equation}
\label{sys1}
\dot{x} = Ax + Bu, \hspace{6mm} y = Cx
\end{equation} 
having $m$ inputs, $p$ outputs, and $n$ states.  The (full)
feedback group is the group generated through the feedback
action:
\begin{equation} \label{1} u\longmapsto u+Fy
\end{equation} 
and through the change of basis in state space, input space and
output space respectively, i.e. through the transformations:
\begin{eqnarray} 
x &\longmapsto & Sx,\hspace{9mm} S\in Gl_n  \label{2} \\ 
u &\longmapsto & T_1u,\hspace{9mm} T_1\in Gl_m  \label{3}\\ 
y &\longmapsto & T_2y,\hspace{9mm} T_2\in Gl_p \label{4}.
\end{eqnarray}

The orbits under the full feedback group are referred to as the
{\em output feedback orbits}. In order to fully understand the effect
of output feedback on the structure of linear systems it is of
fundamental interest to (i) classify the feedback orbits, (ii) to
determine a complete set invariants for output feedback and (iii)
to obtain a detailed description of the adherence order (orbit
closure inclusion) of the different orbits.

Those obviously important problems have already been studied by
many authors (see e.g.~\cite{by85,by86,fu88a,he89a1,hi90,hi92,ra97p})
and despite many partial results the problem is still far from
being solved.

The transformations induced by the
actions\eqr{1},\eqr{2},\eqr{3},\eqr{4} describe a group action on
the vector space of all matrix triples $(A,B,C)$ which is a
vector space of dimension $n(m+n+p)$. There is an extensive
mathematical literature on the classification of orbits arising
from group actions on vector spaces and more general algebraic
varieties and we refer to Section 3 for some more details. If the
number of orbits is finite then this study generally seeks a
discrete set of invariants classifying the finitely  many 
orbits. There are a few instances in the control literature where
the set of orbits is finite and as examples we refer
to~\cite{br70a1,hi92,hi95p}.

In the problem at hand the number of feedback orbits is in
general infinite and this makes the problem difficult. In order
to classify all orbits it will therefore be necessary to derive a
`continuous set of invariants'.

The application of tools from geometric invariant theory (see
e.g.~\cite{mu82,ne78}) often enables one to derive for a given
group action a set of invariants in a systematic way. From a
geometric point of view this amounts to describing an algebraic
variety whose points parameterize uniquely the closed feedback
orbits.

In this paper we construct, using tools from geometric invariant
theory, such a quasi-projective algebraic variety, whose points
parameterize closed output feedback orbits in a unique way.  Since a
quasi-projective variety can be embedded into affine space using
e.g. semi-algebraic functions our result implies the existence of
a complete set of semi-algebraic invariants for output feedback.
It also helps to construct such a complete set of invariants;
however this problem will not be addressed here.

In order to achieve the result it is crucial to first extend the
output feedback action to an action that operates on a
compactification of the space of proper transfer functions of
McMillan degree $n$. This process will be explained in
Section~\ref{Sec:Cas}. In Section~\ref{Sec:Bas} we summarize some
important notions from geometric invariant theory to the extend
we will need it in this paper. 

In order to apply the theorems from geometric invariant theory to
the output feedback invariant problem it will be necessary to
compactify the manifold of $p\times m$ transfer functions of
McMillan degree $n$. This will be accomplished in
Section~\ref{Sec:Pro} using the so called space of homogeneous
autoregressive systems~\cite{ra94}. 

The main results of the paper are provided in
Section~\ref{Sec:main}, where we show that the space of homogeneous
autoregressive systems contains a non-empty Zariski open subset
of semi-stable orbits. This in turn will then lead to a
quasi-projective variety which parameterizes the set of output
feedback invariants in a continuous manner.

In Section~\ref{Sec:Fee} we reinterpret the obtained results in
terms of generalized first order representations. Finally
in~Section~\ref{Sec:Con} we concretely describe the
quasi-projective variety derived in~Section~\ref{Sec:main} in the
situation of single output systems.

\Section{Cascade equivalence and the extended feedback group}
        \label{Sec:Cas}

The notion of {\em cascade equivalence} was introduced by Byrnes
and Helton in~\cite{by86} and it is closely related to the
feedback classification problem.  In our context this notion can
be equivalently described in the following way:

Consider a time invariant linear {\bf proper} system of the form
\begin{equation}
\label{sys2}
\dot{x} = Ax + Bu, \hspace{6mm} y = Cx +Du.
\end{equation} 
In addition to the feedback action\eqr{1} and the basis
transformations\eqr{2},\eqr{3} and\eqr{4} we will also allow a
feed-forward transformation
\begin{equation}   \label{5} 
 y\longmapsto y+Gu.
\end{equation} 
The collection of all those transformations will be called the
{\em extended full feedback group}.  The actions
in\eqr{1},\eqr{3},\eqr{4} and\eqr{5} describe invertible
transformations on the space of external variables $\left[ u^t\ 
  y^t\right]^t$. In this way we can view these actions as
elements of the general linear group $T\in Gl_{m+p}$ and the
collection of these transformations is compactly described
through:
\begin{equation}                \label{6}
\zwei{u}{y}\longmapsto \vier{T_1}{F}{G}{T_2}\zwei{u}{y}
 = T\zwei{u}{y}, 
\end{equation} 
where $ T\in Gl_{m+p}$.  This shows that\eqr{1},\eqr{3},\eqr{4}
and\eqr{5} generate the whole general linear group $T\in
Gl_{m+p}$.

Note that the linear transformation $T\in Gl_{m+p}$ induces the
notion of cascade equivalence on the set of proper transfer
functions. The following Lemma is easily established:
\begin{lem}                          \label{lem:bij}
  There is a bijective correspondence between the set of
  equivalence classes of the form\eqr{sys2} under the extended
  full feedback group and the set of equivalence classes of the
  form\eqr{sys1} under the full feedback group.
\end{lem}

This Lemma now enables us to concentrate on the linear
transformation\eqr{6}.  Instead of working with a state space
description we  can also work with polynomial matrices. For this
let
$$
D^{-1}(s)N(s):= G(s):=C(sI-A)^{-1}B+D
$$
be a left coprime factorization of the transfer function of
system\eqr{sys2}.  Then the linear transformation~\eqr{6} is
equivalently described through:
\begin{equation}                \label{7}
\left(D(s)\ N(s)\right) \longmapsto
 \left(D(s)\ N(s)\right) T^{-1} \hspace{4mm}
T\in Gl_{m+p}.
\end{equation}

\Section{Basic notions from geometric invariant theory}
        \label{Sec:Bas}

Geometric invariant theory constitutes an active research area of
algebraic geometry. One of the main references is the book by
Mumford and Fogarty~\cite{mu82}. The non-specialists among the
interested readers will find the book by Newstead~\cite{ne78} a
good introductory book.

In this section we explain an important result from geometric
invariant theory which we will use later in the paper to derive a
set of continuous feedback invariants. 

Let $X$ be a projective variety, i.e. $X$ is the zero locus of a
finite set of homogeneous polynomial equations. We will assume
that $X$ is embedded into the projective space $\PN$. Let
$G\subset Gl_{N+1}$ be a reductive group (such as, e.g., a group isomorphic
to the general linear group) which acts on the projective space
$\PN$ and induces an action on the variety $X$.

In this situation one has the following general result:
(see~\cite[Theorem~3.14]{ne78}).
\begin{thm}               \label{git}
  There exists a Zariski open set $X^{ss}$ of so called {\em
    semi-stable} points, a projective variety $Y$ and an
  algebraic morphism $\phi:X^{ss}\rightarrow Y$ having the
  property that $\phi^{-1}(y)$ contains exactly one closed $G$
  orbit for every $y\in Y$. Moreover there is a Zariski open set
  $Y^s\subset Y$ such that $\phi^{-1}(y)$ contains one and only
  one orbit for every $y\in Y^s$.
\end{thm} 
The set $X^s:=\phi^{-1}(Y^s)$ is the so called set of stable
points and both $X^s$ and $X^{ss}$ are Zariski open sets of the
variety $X$.  It is possible that $X^{ss}$ is the empty set in
which case Theorem~\ref{git} does not give any insight.

Theorem~\ref{git} is significant in several ways. First the
points of the variety $Y$ provide a continuous family of
invariants capable of distinguishing orbits inside $X^s$. In
addition the variety $Y$ is characterized through some universal
properties. Because of this reason one sometimes also speaks
about the {\em categorical quotient} $Y$.

The fact that this variety $Y$ is projective is surprising.  It
will be our goal in the next section to apply Theorem~\ref{git}
to the feedback orbit classification problem.

\Section{The projective variety of homogeneous autoregressive
  systems}   \label{Sec:Pro}

In order to apply Theorem~\ref{git} it will therefore be necessary to
compactify the space of all $p\times (m+p)$ autoregressive
systems of the form $P(s)=\left(D(s)\ N(s)\right)$. Such a
compactification was provided in~\cite{ra94} and we shortly
review the details.

Consider a $p\times (m+p)$ polynomial matrix
\begin{equation} 
P(s,t):= \left( \begin{array}{ccc}
f_{1,1}(s,t)&\ldots&f_{1,(m+p)}(s,t)\\ 
f_{2,1}(s,t)&\ldots&f_{2,(m+p)}(s,t)\\
\vdots&&\vdots\\ f_{p,1}(s,t) &\ldots&f_{p,(m+p)}(s,t)
\end{array}\right).
\end{equation}

We say $P(s,t)$ is homogeneous of row degrees
$ \nu_1, \dots, \nu_p$ if each element $f_{i,j}(s,t)$ is
a homogeneous polynomial of degree $\nu_i$. 
A square matrix $U(s,t)$ of homogeneous polynomials is called
unimodular, if $\det U(s,t)$ is a nonzero monomial in $t$.
 We say two
homogeneous matrices $P(s,t)$ and $\tilde{P}(s,t)$ are equivalent
if they have the same row-degrees and if there is a unimodular
matrix $U(s,t)$, whose entries $u_{ij}(s,t)$ are homogeneous
polynomials of degree $\nu_i-\nu_j$ with $\tilde{P}(s,t) = U(s,t)P(s,t)$.
 Using this equivalence relation we define:
\begin{de}
  An equivalence class of full rank homogeneous polynomial
  matrices $P(s,t)$ will be called a {\em homogeneous
    autoregressive system}. The McMillan degree of a homogeneous
  autoregressive system is defined as the sum of the row degrees,
  i.e. through $n:=\sum_{i=1}^{p}\nu_i$.  The set of all
  homogeneous autoregressive systems of size $p\times (m+p)$ and
  McMillan degree $n$ will be denoted by ${\mathcal H}^n_{p,m}$.
\end{de}

Let ${\mathcal R}_{n,m,p}$ denote the space of $p\times m$ proper
transfer functions of McMillan degree $n$.  The main result
established in~\cite{ra94} is as follows:
\begin{thm}[\cite{ra94}]  \label{main}
  ${\mathcal H}^n_{p,m}$ is a smooth projective variety
  containing the set of proper transfer functions
  ${\mathcal R}_{n,m,p}$ as a Zariski dense subset.
\end{thm}

More generally, ${\mathcal H}^n_{p,m}$  contains the set of all
degree $n$ rational curves of the Grassmannian Grass $(m,m+p)$ as
a Zariski-dense subset.
The variety ${\mathcal H}^n_{p,m}$ arises in algebraic geometry
in the following context: Let $\pp^1$ be the projective line and
let $O_{\pp^1}$ be the structure sheaf of $\pp^1$. Let $V$ be an
$(m+p)$-dimensional vector space over the base field $K$.  The
space ${\mathcal H}^n_{p,m}$ is the quotient scheme that
parameterizes all quotients $\mathcal B$ of the sheaf $V\otimes
O_{\pp^1}$ of degree of  $n$ and 
rank  $m$, that is, sheaves $\mathcal B$ whose Hilbert
polynomial is $\chi(\ell)=m(\ell+1)+n$. This identification
proceeds as follows(refer to \cite{ra94} for more details): A
point $x$ in the quotient scheme gives rise to a short exact
sequence:
\begin{equation}\label{ses} 
0\to{\mathcal A}\overset{\psi}{\to} V\otimes
O_{\pp^1}\overset{\phi}{\to} {\mathcal B}\to 0
\end{equation}

Now $\mathcal A$ is a locally free sheaf of degree $-n$ and rank
$p$. So ${\mathcal A}\simeq \underset{i=1}{\overset{l}\bigoplus}
O_{\pp^1}(-\nu_i)$.  Therefore the map from $\mathcal A$ to
$V\otimes O_{\pp^1}$ is given by the transpose of a homogeneous
autoregressive system of the form $P(s,t)$, once a basis is
chosen for $V$ and $\mathcal A$. Conversely given a homogeneous
autoregressive system $P(s,t)$ the map defined by the transpose
of $P(s,t)$ from $\underset{i=1}{\overset{l}\bigoplus}
O_{\pp^1}(-\nu_i)$ to $V\otimes O_{\pp^1}$ is an injective map
since $P$ has full rank. The quotient of this map has rank $m$
and degree $n$ and thus defines a point $x$ in the quotient
scheme.

The space ${\mathcal H}^n_{p,m}$ can be embedded as a projective
variety in the following way: Fix an integer $\ell\ge n$. Given a
point $x\in {\mathcal H}^n_{p,m}$ corresponding to the short
exact sequence of sheaves~(\ref{ses}) one has the corresponding
long exact sequence of cohomology groups:
$$
0\to H^0(\pp^1, {\mathcal A}(\ell))\to H^0(\pp^1, V\otimes
O_{\pp^1}(\ell))\to H^0(\pp^1, {\mathcal B}(\ell))\to H^1(\pp^1,
{\mathcal A}(\ell))\to \cdots
$$
Now, if $\ell\ge n$ then $H^1(\pp^1, {\mathcal A}(\ell))=0$.
Also $\dim H^0(\pp^1, {\mathcal A}(\ell))=p(\ell+1)-n$ and
$H^0(\pp^1, V\otimes O_{\pp^1}(\ell))=(m+p)(\ell+1)\simeq
V\otimes H^0(\pp^1,O_{\pp^1}(\ell))$, so that 
$\dim H^0(\pp^1, V\otimes O_{\pp^1}(\ell)) =(m+p)(\ell +1)$.

  Finally one obtains a map
$\rho_{\ell}: {\mathcal H}^n_{p,m}\to \mbox{Grass} (p(\ell+1)-n,
V\otimes H^0(\pp^1,O_{\pp^1}(\ell))$, the Grassmanian of
$p(\ell+1)-n)$-dimensional subspaces of  $V\otimes
H^0(\pp^1,O_{\pp^1}(\ell)$,
obtained by defining
$\rho_{\ell}(x)$ to be the subspace $H^0(\pp^1, {\mathcal
  A}(\ell))$. The map $\rho_{\ell}$ defines an embedding (see
\cite{ra94} where it is proved specifically for $\ell=n$, but the
same proof applies to $\ell>n$). This Grassmannian can be
embedded through the Pl\"ucker embedding in 
$\pp= \pp(\overset{(p(\ell+1)-n)}\wedge V \otimes
H^0(\pp^1,O_{\pp^1}(\ell)))$, the projective space of lines in this vector
space.
 The group $\mbox{Gl}(V)$
obviously acts on $V$, therefore also on $V\otimes O_{\pp^1}$ and
thus also on the vector space $V\otimes
H^0(\pp^1,O_{\pp^1}(\ell))$ for each $\ell$ and also on
$\overset{(p(\ell+1)-n)}\wedge V\otimes
H^0(\pp^1,O_{\pp^1}(\ell))$. Thus for each $\ell\ge n$ there is
an induced action of $\mbox{Gl}(V)$ on ${\mathcal H}^n_{p,m}$ as
an embedded subvariety of the projective space $\pp$.

\Section{Main Results}    \label{Sec:main}

Our first step will be to identify the semi-stable points for
this action of $\mbox{Gl}(V)$. The main technical tool we will
use throughout this section is the following result of
Simpson(\cite{si94}):

\begin{thm}[Simpson, \cite{si94}, Lemma 1.15] \label{crit} 
  There exists an $L$ such that for $\ell\ge L$ the following
  holds: Suppose $x: V\otimes O_{\pp^1}\to{\mathcal B}\to 0$ is a
  point in ${\mathcal H}^n_{p,m}$. For any subspace $H\subset V$,
  let $\mathcal G$ denote the subsheaf of $\mathcal B$ generated
  by $H\otimes O_{\pp^1}$. Suppose that $\chi ({\mathcal
    G},\ell)>0$ and 
  \begin{equation}\label{simpss}
    \frac{\mbox{dim}\,H}{\chi({\mathcal
        G},\ell)}\le\frac{\mbox{dim}\,
      V}{(m(\ell+1)+n)}
  \end{equation} 
  (resp. $<$) for all nonzero proper subspaces $H \subset V$.
  Then the point $x$ is semi-stable (resp. stable) in the
  embedding $\rho_\ell: {\mathcal H}^n_{p,m}\to \pp^{N}$
  described above.
\end{thm}  

We want to rephrase this theorem in a geometric form that will be
more suitable for our application. As a first step we want to
point out that every point in the compactification may be viewed
as a map from the projective line $\pp^1$ to a Grassmannian.
While this identification is fairly standard for the points $x\in
{\mathcal H}^n_{p,m}$ where the quotient sheaf $\mathcal B$ is
locally free we want to explain how it can be extended to the
points where $\mathcal B$ is not locally free.
  
If $x\in {\mathcal H}^n_{p,m}$ corresponds to the short exact
sequence~(\ref{ses}) then ${\mathcal B}\simeq {\mathcal B}_{\mbox
  {free}} \oplus{\mathcal
  B}_{\mbox {tor}}.$ %
The map $\phi$ followed by the surjection: %
${\mathcal B}_{\mbox {free}}\oplus{\mathcal B}_{\mbox {tor}}\to
{\mathcal B}_{\mbox {free}}\to 0$ %
gives rise to a surjection: %
$x':\, V\otimes O_{\pp^1}\overset{\phi'}{\to}{\mathcal
  B}_{\mbox{free}}$ %
which corresponds to an observable system of the same rank, but
lower degree than the original system $x$. The point $x'$ belongs
to ${\mathcal H}^{n'}_{p,m}$ where $n'<n$. The system $x'$ is the
observable part of the system $x$.

Let $G'=\mbox{Grass}(m,V)$ be the Grassmannian of $m$-dimensional

{\em quotients} of $V$ then each point $x\in{\mathcal H}^n_{p,m}$
corresponds to a map $\phi_x: \pp^1\to G'$ given as follows: If
$x$ corresponds to the short exact sequence: $0\to{\mathcal A}\to
V\otimes O_{\pp^1}\overset{\phi}\to{\mathcal B}\to 0$, then we
have a surjection $\phi$ (or $\phi'$) from $V\otimes O_{\pp^1}$
to $\mathcal{B}_{\mbox {free}}$. This map can be represented by a
$m\times (m+p)$ matrix $Q(s,t)$ where each row of this matrix
consists of homogeneous polynomials $g(s,t)$ of the same degree.
Further at each point $(s,t)\in \pp^1$ the matrix $Q(s,t)$ has
full row rank. Thus given a point $z\in \pp^1$ we can define a
quotient denoted by $\phi(z)$ of rank $m$ of the vector space $V$
given by the matrix $Q(s,t)$ evaluated at the point $z$. We
denote the corresponding map from $\pp^1$ to $G'$ by $\phi$.

\begin{lem}\label{genss}
  Let $x\in{\mathcal H}^n_{p,m}$. Suppose that for a generic
  point $z\in\pp^1$ and any proper subspace $H\subset V$
  \begin{equation}\label{ss}
  \dim (\phi_x(z))(H)> \frac{m}{m+p}\dim H
  \end{equation}
  where $\phi_x$ is the map from $\pp^1$ to the Grassmannian $G'$
  of quotients of $V$, associated to the point $x$ and $\phi_x(z)(H)$
denotes the image of $H$ under the canonical projection map $V
  \rightarrow \phi_x(z)$. Then there
  exists an $L$ such that for $\ell\ge L$, $x$ is a stable point
  in $\rho_{\ell}({\mathcal H}^n_{p,m})$.
\end{lem}

\begin{proof}
  In the first part of the proof we will assume that the point
  $x$ is observable, that is, the sheaf $\mathcal B$ is locally
  free.
  
  Let $H\subset V$ be a proper subspace of $V$ and let the image
  of $H$ under the map $\phi$ be $\mathcal{G}$. If~(\ref{ss}) is
  satisfied then the rank $g$ of the sheaf $\mathcal G$, which is
  equal to the dimension of $(\phi_x(z))(H)$ at the generic point
  is greater than $\frac{m}{m+p}\dim H$. Now the Euler
  characteristic $\chi ({\mathcal G}(\ell))=g(\ell+1)+\deg
  {\mathcal G}$. So for $\ell$ large enough
  $$
  \begin{array}{ccc}
  \chi ({\mathcal G}(\ell))&> &\dim H\frac{m}{m+p}(\ell+1)+n\\
   &>& \frac{\dim H}{m+p}(m(\ell +1)+n)\\
  \mbox{ so } \frac{\dim H}{\chi ({\mathcal G}(\ell))}&<&
  \frac{m+p}{m(\ell +1)+n)}.
  \end{array}
  $$
  Therefore by Theorem~\ref{crit} the point $x$ is stable.
  
  If $x$ is not observable our assumption is that the map
  $\phi'_{x'}$ satisfies condition~(\ref{ss}). The rank of
  $\phi'_{x'}$ is the rank of the sub-sheaf ${\mathcal G}'$ of
  ${\mathcal B}_{\mbox {free}}$ generated by $H\otimes
  O_{\pp^1}$.  The sheaf ${\mathcal G}'$ is a sub-sheaf of the
  sheaf $\mathcal G$ which is the sub-sheaf of $\mathcal B$
  generated by $H\otimes O_{\pp^1}$, so $\chi({\mathcal G}
  (\ell))\ge \chi({\mathcal G}' (\ell))$. By the above
  calculations $\chi({\mathcal G'} (\ell))$ satisfies
  condition~\ref{simpss} so $\mathcal G$ satisfies this condition
  as well. Thus the point $x$ is stable.
\end{proof}

We wish to recall the following definition from the systems
theory literature:
\begin{de}
  A $p\times (m+p)$ homogeneous autoregressive system $P(s,t)$ is
  called {\em nondegenerate} if there is no full rank $m\times
  (m+p)$ matrix $K$ with entries in the complex numbers $\C$ such
  that
  $$
  \det \genfrac{(}{)}{0pt}{0}{P(s,t)}{K}=0.
  $$
\end{de}

One verifies that nondegenerate systems cannot exist if the
McMillan degree is `small'. The following result shows when
nondegenerate systems are open and dense inside the variety
${\mathcal H}^n_{p,m}$.
\begin{thm}[\cite{br81}]
  If the McMillan degree $n\geq mp$ then the variety ${\mathcal
    H}^n_{p,m}$ contains a nonempty Zariski open set of nondegenerate
  systems.
\end{thm}

\begin{lem}      \label{lem:sta}
  If $P(s,t)$ is a homogeneous autoregressive system that is
  nondegenerate then the corresponding point $x\in{\mathcal
    H}^n_{p,m}$ satisfies the condition~(\ref{ss}) of
  Lemma~\ref{genss} and is therefore a stable point in
  $\rho_{\ell}({\mathcal H}^n_{p,m})$ for large enough $\ell$.
\end{lem}

\begin{proof} 
  Suppose the point $x\in {\mathcal H}^n_{p,m} $ corresponds to
  the short exact sequence~(\ref{ses}). Then the map $\psi$ is
  represented by the transpose of the matrix $P(s,t)$. For each
  point $z\in \pp^1$, $\psi(z)$ determines a subspace of $V$
  given by the row span of the matrix $P(s,t)$ evaluated at the
  point $z$. If $P(s,t)$ is nondegenerate then for the generic
  point $z\in \pp^1$ and a subspace $H\subset V,\, \dim
  (H\cap\psi (z))=0$ if $\dim H\le m$ and if $\dim H>m$ then
  $\dim (H\cap\psi (z))=\dim H-m$. Notice that $\psi(z)$ is the
  kernel of the map $\phi(z)$. So if $x$ is a nondegenerate point
  at the generic point $z\in \pp^1$, if $\dim H\le m$ then $\dim
  \phi_x(H)=\dim H$ and if $\dim H>m$ then $\dim (\phi_x(H))=m$.
  Thus the point $x$ satisfies the condition~(\ref{ss}) and
  therefore it is stable.
\end{proof}

We want to remark that the converse of the statement in the above
Lemma is not true. One can find stable points in
$\mathcal{H}^n_{p,m} $ that are not nondegenerate 
as the following example shows.
\begin{exmp}                    \label{degst} 
  Let $P(s,t)$ be given by the following matrix:
\begin{equation}\label{ex}
P(s,t)=\left(\begin{array}{ccccc}
s^2&st&t^2&s^2&s^2+t^2\\
st&t^2&s^2&s^2+2t^2&s^2-t^2
\end{array}\right)
\end{equation}
where the transpose of $P(s,t)$ is the matrix of the sheaf map
$\psi$ in~(\ref{ses}).  Let $x\in {\mathcal H}^4_{3,2}$ be the
point represented by the homogeneous autoregeressive system
$P(s,t)$. This point is degenerate since
$$
\det \left(\begin{array}{c}
P(s,t)\\
{\begin{array}{ccccc}0&0&1&0&0\\0&0&0&1&0\\
0&0&0&0&1\end{array}}
\end{array}\right)=0.
$$
 The matrix of the map $\phi (s,t)$ in the short exact
 sequence~(\ref{ses}) corresponding to this point $x$
 was computed to be the following matrix:
 $$Q(s,t)=\left(\begin{array}{ccccc}
 -t&s&0&0&0\\
 s+5t&-2t&s+4t&-2s-t&s-4t\\
 -s^2-4st-8t^2&3t^2&-s^2-3st-7t^2&s^2+4st+2t^2&7t^2\end{array}\right).
$$
In order to check that the point $x$ is stable, according to
Lemma~\ref{genss} it suffices to check that for a generic point
$z\in \pp^1$, for any subspace $H$ of dimension three the image
$\phi_x(z)(H)$ has dimension at least two, and if $\dim H$ is
four then the image has dimension at least three.  Let $A\in
Gl_5$ be a generic matrix. Now, it is enough to check that the
first three columns of the matrix $QA$ has rank at least two, and
the first four columns has a rank of at least three, at the
generic point in $\pp^1$. We confirmed this to be the case using
the computer algebra system Maculay 2, even though the second
computation took more than an hour on a Macintosh G4.  Thus the
point $x$ is degenerate but stable.
\end{exmp}

An immediate consequence of Lemma~\ref{lem:sta} is the following
theorem:

\begin{thm}                       \label{thm:sem}
  The set of semi-stable orbits $X^{ss}\subset{\mathcal
    H}^n_{p,m}$ contains the set of nondegenerate systems. In
  particular if $n\geq mp$ then $X^{ss}$ is nonempty and there
  exists a projective variety $Y$ and a morphism
  $\phi:X^{ss}\rightarrow Y$ having the property that
  $\phi^{-1}(y)$ contains exactly one closed $Gl_{m+p}$ orbit for
  every $y\in Y$.
\end{thm}
\begin{proof}
 The stable orbits (and nondegenerate systems are stable by
 Lemma~\ref{lem:sta}) form a subset of the semi-stable orbits.
\end{proof}
We want to remark that since the stabilizer group of a stable
point is a finite group (see~\cite{ne78}) the stabilizer group of
a nondegenerate point is finite.

\begin{cor}                        \label{cor:sem}
  Consider the set of $m$ input, $p$ output systems of McMillan
  degree $n$ as introduced in\eqr{sys1}. If $n\geq mp$ then the
  set of feedback orbits with respect to the full feedback group
  has a continuous set of invariants consisting of a
  quasi-projective variety.
\end{cor}
\begin{proof}
  Consider the set ${\mathcal R}_{n,m,p}$ of proper $p\times m$

  transfer functions of McMillan degree $n$. Let
  $V:= {\mathcal R}_{n,m,p}\cap X^{ss}$ and let $\hat{\phi}$ be the
  restriction of the morphism $\phi$ to the Zariski open subset
  $V$. Then $\mathrm{im}(\hat{\phi})$ describes a
  quasi-projective variety parameterizing the orbits under the
  extended full feedback group as described in
  Section~\ref{Sec:Cas}. By Lemma~\ref{lem:bij} this variety also
  parameterizes the set of $m$ input, $p$ output systems of
  McMillan degree $n$ as introduced in\eqr{sys1} modulo the full
  feedback group.
\end{proof}

\Section{The feedback orbit classification problem in terms of
  generalized first order systems}     \label{Sec:Fee}

The set of homogeneous autoregressive systems can be described
through  generalized first order systems and we refer
to~\cite{ra94,ra97}. In this section we describe the extended
feedback group in terms of these generalized first order system.
In this way we make the connection with work of Hinrichsen and
O'Halloran~\cite{hi95p}.

Following the exposition in~\cite{ku94,ra95} consider
$(n+p)\times n$ matrices $K,L$ and a $(n+p)\times (m+p)$ matrix
$M$.  Those matrices define a generalized state space system
through
\begin{equation}            \label{dualsystem}
 K\dot{x}(t)+Lx(t)+Mw(t)=0,
\hspace{5mm} x(t)\in \R^{n},w(t)\in \R^{m+p}.
\end{equation}

The system\eqr{dualsystem} is called admissible if the
homogeneous pencil $\left[ {sK\, +\, tL}\right]$ has generically
full column rank.  An admissible system is called controllable if
the pencil $\left[{sK+ tL\; M}\right]$ has full row rank for all
$(s,t) \in \C^2 \setminus \{(0,0)\}$.

There is a natural equivalence relation among generalized first
order systems: If $U\in Gl_{n+p}$ and $S\in Gl_{n}$ then
\begin{equation}               \label{sim2} (K,L,M)
\hspace{2mm}\sim\hspace{2mm}(UKS^{-1},ULS^{-1},UM) .
\end{equation} 
If the high order coefficient matrix
$$
[K\, M]
$$
has the property that the first $(n+p)\times (n+p)$ minor is
invertible then the system\eqr{dualsystem} is equivalent to a
system having
\begin{equation}             \label{ABCD2} 
K=\left[ \begin{array}{c}-I\\ 0
\end{array}\right], \hspace{3mm} 
L=\left[ \begin{array}{c}A\\ C
\end{array}\right], \hspace{3mm} 
M=\left[ \begin{array}{cc}0&B\\ -I&D
\end{array}\right],
\end{equation} 
i.e. the system is equivalent to a usual state space system of
the form
$$
\dot{x}=Ax+Bu,\; y=Cx+Du.
$$

The connection to the set of homogeneous autoregressive
systems is established through:

\begin{thm} (\cite{ra95})      \label{main-co}
  The categorical quotient of the set of controllable state space
  systems as introduced in\eqr{dualsystem} under the group action
  $Gl_{n+p}\times Gl_n$ is isomorphic to the smooth projective
  variety ${\mathcal H}^{n}_{m,p}$ of all $m\times (m+p)$
  homogeneous autoregressive systems of McMillan degree $n$.
\end{thm}

Next we define:
\begin{de}
  Two generalized first order systems are equivalent under the
  extended full feedback group if there are invertible matrices
  $S\in Gl_{n}$, $T\in Gl_{m+p}$ and $U\in Gl_{n+p}$ such that
\begin{equation}    \label{10} 
(K,L,M) \hspace{2mm}\sim\hspace{2mm}(UKS^{-1},ULS^{-1},UMT^{-1}) .
\end{equation}
\end{de}

Note that the linear transformation $T$ introduced in~\eqr{6}
corresponds to a change of basis in the set of external variables
$w=\zwei{u}{y}$ and it is therefore equal to the transformation
$T$ appearing in\eqr{10}. Finally the group action described
in\eqr{10} corresponds exactly to the transformations
$(i),(ii),(iii)$ considered by Hinrichsen and O'Halloran
in~\cite[p. 2730]{hi95p}.

As a consequence of Theorem~\ref{main-co} we have:

\begin{thm}
  The categorical quotient induced by the group action\eqr{10} is
  equal to the projective variety ${\mathcal H}^{n}_{m,p}/
  Gl_{m+p}.$ Moreover ${\mathcal H}^{n}_{m,p}/ Gl_{m+p}$
  represents a continuous parameterization of the feedback orbits
  under the full feedback group.
\end{thm}

\Section{A concrete description of the moduli space in the MISO
  situation}  \label{Sec:Con}

In this section we explain our result in the multi input, single
output (i.e.  $p=1$) situation. The variety ${\mathcal
  H}^{n}_{1,m}$ consists in this case of all $1\times (m+1)$
polynomial vectors
$$
P(s)=(f_1(s),\ldots,f_{m+1}(s))
$$
whose polynomial entries have degree at most $n$. In this way
we can identify the space ${\mathcal H}^{n}_{1,m}$ with the
projective space $\PN$, where $N=mn+m+n$.

By definition a systems is nondegenerate if the $m+1$ polynomial
vectors 
$$
\{ f_1(s),\ldots,f_{m+1}(s)\}\subset\R[s]
$$
are linearly independent over $\R$.  Clearly this can only
happen if the McMillan degree $n\geq m$. Theorem~\ref{thm:sem}
and Corollary~\ref{cor:sem} state in this case that the set of
systems where $\{ f_1(s),\ldots,f_{m+1}(s)\}$ are linearly
independent are all contained in the semi-stable orbits and that
there is a quasi-projective variety describing the quotient. In
our situation this can be made very concrete:

Identify the set of polynomial vectors of degree at most $n$ with
the vector space $\R^{n+1}$. A system
$P(s)=(f_1(s),\ldots,f_{m+1}(s))$  then defines a linear subspace
$$
{\rm span}_\R \{ f_1(s),\ldots,f_{m+1}(s)\}\subset\R^{n+1}.
$$
This subspace has dimension $m+1$ if and only if $P(s)$
describes a nondegenerate system. The semi-stable points under the
extended full feedback group $Gl_{m+1}$ therefore describe  a well
defined $m+1$ dimensional subspace of $\R^{n+1}$. The categorical
quotient ${\mathcal H}^{n}_{1,m}/ Gl_{m+1}$ is in this case
exactly the Grassmann variety Grass$(m+1,\R^{n+1})$ of $m+1$
dimensional subspaces in $\R^{n+1}$. This extends a construction in
[3,6] for output feedback invariants of SISO systems, where $m = 1.$

In particular it follows that the semi-stable orbits coincide
exactly with the set of nondegenerate systems, something which is
not true for general multi-output systems, as illustrated by 
Example~\ref{degst}.

The case $p = 0$ is interesting and nontrivial as well. In fact, the
smooth moduli space, as constructed above, parametrizes
representations for the wild input-state quiver action $(K,L,M)
\mapsto (U K S^{-1}, ULS^{-1}, UMT^{-1})$. The task of classifying
representations of this quiver has been an open problem for at least
the last two decades.

\Section{Conclusion}

In this paper we did show that the output feedback problem is
closely related to the study of the moduli space
$\mathcal{H}^{n}_{p,m}/ Gl_{m+p},$ where $\mathcal{H}^{n}_{p,m}$
denotes the set of $p\times (m+p)$ homogeneous autoregressive
systems of McMillan degree $n$. 

\nocite{ku94,le97b,wo76,mo73}
\end{document}